\documentclass[a4paper,leqno]{article} \usepackage{amsmath,amssymb,small}

\author{David R. E. Williams\footnote{The author was supported by EPSRC grant GRL67899.}}
\title{Diffeomorphic flows driven by\\ L\'evy processes} 
\date{}

\begin{document}
\maketitle
\def\thefootnote{\fnsymbol{footnote}}
\def\@makefnmark{\hbox to\z@{$\m@th^{\@thefnmark}$\hss}}
\footnotesize\rm\noindent
\hspace*{16pt}\strut\footnote[0]{{\it AMS\/\ {\rm 1991} subject
classifications}.  60H20}
\normalsize\rm

{\abstract We prove that the stochastic differential equation 
\begin{equation}\label{e.0}
Y_{s,t}(x) = Y_{s,s}(x) + \int_0^{t-s} f(Y_{s,s+u}(x))\;dX_{s+u}, \qquad Y_{s,s}(x)=x\in\R^d.
\end{equation}
driven by a L\'evy process whose paths have finite $p$-variation almost surely for some $p\in[1,2)$ defines a flow of locally $C^1$-diffeomorphisms provided the vector field $f$ is $\alpha$-Lipschitz for some $\alpha>p$.  Using a path-wise approach we relax the smoothness condition normally required for a class of discontinuous semi-martingales.}
\section*{Introduction}\label{s.intro}
In this paper we give sufficient conditions for the solution  of the following integral equation to define a flow of local $C^1$-diffeomorphisms on $\R^d$ when the integrator is allowed to have unbounded variation:
\begin{equation}\label{e.1}
Y_{s,t}(x) = Y_{s,s}(x) + \int_0^{t-s} f(Y_{s,s+u}(x))\;dX_{s+u}, \qquad Y_{s,s}(x)=x\in\R^d.
\end{equation}

The above assumption means that one must take care defining the integral.  In this paper we consider integrators with finite $p$-variation for some $p\in[1,2)$.  We use the Young \cite{you} integral
\begin{equation}
\int_0^t f(u)dg(u)\qquad\qquad \text{(Young),}
\end{equation}
which is defined by taking a limit of Riemann sums whenever the functions $f$ and $g$ have no common discontinuities and have finite $p$ and $q$ variation for some $1/p+1/q>1$.  Lyons \cite{tel1} solved the following integral equation:
\begin{equation}\label{e.terry}
y_t = y_0 + \int_0^t f(y_u)\,dx_u\qquad y_0=a\in\R^d,
\end{equation}  
where $x_t$ is a continuous function, having finite $p$-variation for some $p<2$ and $f$ is a vector field which is $\alpha$-Lipschitz for some $\alpha>p-1$.  Furthermore, he proved that if $\alpha>p$ then the solution is unique.

In this paper we show that the solutions of equation \eqref{e.terry} form a flow of local $C^1$-diffeomorphisms.  Then we extend the solution of equation \eqref{e.terry} to allow the integrator to be discontinuous.  Finally we apply the results to SDEs driven by a class of L\'evy processes whose paths, almost surely, have finite $p$-variation for some $p<2$.  It is shown (Corollary \ref{cor.levy}) that one can relax the smoothness condition used when using the It\^o calculus within the semi-martingale framework.

The paper is organised in the following way:  The first section introduces the notation used in the remainder of the paper.  The second section contains the proof of the diffeomorphic flow property.  In the third section we extend the result of Lyons to incorporate integrators with discontinuities.  Finally we apply the results to the class of L\'evy processes whose paths a.s. have finite $p$-variation for some $p<2$.

\section{Preliminaries}\label{s.prelim}

In this section we recall some definitions which will be required throughout the paper.

\begin{defn}
The $p$-variation of a function ${ X }(s)$ over the interval $[0,t]$ is defined as follows:
\begin{equation}
{\Vert}{ X }{\Vert}_{_{p,[0,t]}}=\bigg\{\sup_{\pi\in\pi[0,t]}\quad{\sum_\pi}{\vert}{ X }(t_k)-{ X }(t_{k-1}){\vert}^p\bigg\}^{1\over
p}
\end{equation}
where $\pi[0,t]$ is the collection of all finite partitions of the interval
$[0,t]$.  
\end{defn}
\begin{rems}
Note that ${\Vert}.{\Vert}_p$ is a semi-norm on the class of
functions whose p-variation is finite.  Let {\bf $W_p$} denote this class of
functions.  Trivially one has {\bf $W_p$} $\subseteq$ {\bf $W_q$} when $p\leq q$.  Let $W_p({ a })$ denote the subspace of paths in $W_p$ with initial point ${ a }$.  Then $W_p({ a })$ equipped with $\norm{.}{p}$ is a Banach space.  

Let ${ X }_t$ be a path of finite $p$-variation on $[0,T]$.  Then necessarily
there are only countably many discontinuities and there exist left and
right limits of ${ X }_t$ for all $t\in [0,T]$.

Note that we use a stronger version of $p$-variation than is used in the stochastic calculus literature.  Fristedt \cite{fris} compares the different definitions of variation.
\end{rems}
\begin{eg}
Probabilists are used to the idea that Brownian motion a.s. has finite quadratic variation.  The process a.s. fails to have finite 2-variation in the above (strong) sense; it does however have finite $p$-variation a.s. for all $p>2$.  
\end{eg}
\begin{defn}
  A function $f$ is in $\lip(\alpha)$ for some $\alpha>1$ if
\begin{equation*}
\norm{f}{\infty} < \infty \;\hbox{and}\; {{\partial f}\over{\partial x_j}} \in \hbox{Lip}(\alpha -1)\quad j=1,\dots ,d\,.
\end{equation*}
Its norm is given by
\begin{equation*}
\norm{f}{\lip(\alpha)}  \dfn  \norm{f}{\infty} + \sum_{j=1}^d \lnorm{{{\partial{f}}\over{\partial{x_j}}}}{\lip(\alpha -1)}	\qquad\qquad\hbox{for}\;\alpha>1.
\end{equation*}
\end{defn}

This is Stein's \cite{ste} definition of $\alpha$-Lipschitz continuity for $\alpha>1$.  It extends the classical definition:  $f$ is in $\lip(\alpha)$ for some $\alpha\in(0,1]$ if
\begin{equation*}
\snorm{f(x) - f(y)} \leq K \snorm{x-y}^{\alpha}
\end{equation*}
with norm
\begin{equation*}
\norm{f}{\infty} + \sup_{x \neq y} \frac{\snorm{f(x)-f(y)}}{\snorm{x-y}^{\alpha}}.
\end{equation*}

\section{Diffeomorphic flows}\label{s.flow}

Consider the following differential equation:
\begin{equation}\label{e.ode}
dy_t = \sum_{i=1}^d f^i(y_t)\,dx_t^i \qquad y_0=a\in R^d
\end{equation}
where $(x_t)_{t\geq0}$ is a continuous $\R^n$-valued function with finite $p$-variation for some $p<2$ and $f:\R^d\rightarrow\hom(\R^n,\R^d)$.  Lyons proved the following theorem:
\begin{thm}[\cite{tel1}]
If $f$ is $\alpha$-Lipschitz for some $\alpha>p-1$ then solutions to \eqref{e.ode} exist.  If $f$ is $\alpha$-Lipschitz for some $\alpha>p$ then there is a unique solution.
\end{thm}
\begin{rem}
The existence result proves that the solution path has finite $p\prime$-variation for some $p\prime>p$.  The uniqueness result proves that the solution has finite $p$-variation.  Both results are proved by applying fixed point theorems to a sequence of Picard iterations.
\end{rem} 

We modify the proof of the above theorem to show that a flow of local $C^1$-diffeomorphisms exist.  First, we recall the definition of a flow of diffeomorphisms.

\begin{defn}
A collection of maps $\{\phi_{s,t} \,:\, 0\leq s\leq t \}$ form a flow of homeomorphisms on $\R^d$ if:
\begin{itemize}
\item[(i)] $\phi_{s,u} = \phi_{s,t}\circ\phi_{t,u}$ for all $s\leq t\leq u$.
\item[(ii)] $\phi_{s,s}$ is the identity map on $\R^d$ for all $s\in\R^+$.
\item[(iii)] $\phi_{s,t}$ is an onto homeomorphism of $\R^d$ for all $s,t\in\R^+$.

Further if $\phi_{s,t}$ satisfies (iv) then it is called a flow of $C^k$-diffeomorphisms.
\item[(iv)] $\phi_{s,t}(x)$ is $k$ times differentiable with respect to $x$.  The derivatives are continuous with respect to $(s,t,x)$.
\end{itemize}

\end{defn}
\begin{rem}
Note that a stochastic flow of homeomorphisms \cite[p.114]{kunita} is a continuous random field which satisfies the above properties off a set of measure zero.
\end{rem}

\begin{thm}\label{thm.main}
Let $(X_t)_{t\geq0}$ be a $\R^n$-valued continuous path having finite $p$-\linebreak variation for some $p<2$.  Let $f:\R^d\rightarrow\hom(\R^n,\R^d)$ be $\alpha$-Lipschitz for some $\alpha>p$.  Then the unique solution to the following system of differential equations is a flow of local $C^1$-diffeomorphisms on $\R^d$:
\begin{equation}\label{e.flo}
Y_{s,t}(x) = Y_{s,s}(x) + \int_0^{t-s} f(Y_{s,s+u}(x))\;dX_{s+u}, \qquad Y_{s,s}(x)=x\in\R^d.
\end{equation}
\end{thm}
The proof is an extension of the method used in \cite{tel1}.  The key step is to consider the following pair of differential equations:
\begin{align}\label{e.pair}
dY_{s,u}(x) &= f(Y_{s,s+u}(x))\;dX_{s+u} \qquad Y_{s,s}(x)=x\in\R^d\nonumber\\
dK_{s,u}(x) &= \nabla f(Y_{s,u}(x))\,K_{s,u}(x)\;dX_{s+u} \qquad K_{s,s}(x) = I \in M_{d\times d}
\end{align}
where $M_{d\times d}$ denotes the space of $d\times d$ matrices over $\R$.

\proof
We rewrite \eqref{e.pair} as follows:
\begin{align}\label{e.deriv}
dZ_{s,u}(x) &= h(Z_{s,u}(x))\;dX_{s+u} \qquad Z_{s,s}(x)=(x,I)\in \R^d\times M_{d\times d},
\intertext{where $h$ is the vector field defined by}
h\big((y,k)\big) &\dfn \big(f(y),\nabla f(y) k\big).
\end{align}
Given that $f$ is $\alpha$-Lipschitz for some $\alpha>p$ we know that $h$ is $(\alpha-1)$-Lipschitz.  By applying Lyons' existence result \cite{tel1} we deduce that there exist solutions to \eqref{e.deriv}.  However, the solution of the first equation in \eqref{e.pair} is unique.  This implies that the solution of \eqref{e.deriv} is unique.

It remains to show that the solution is $C^1$.  Let $a,\,b\in\R^d$.  Let $q\in(p,2)$ and $r=(\alpha-1)/p -\delta$ where $\delta\in(0,(\alpha-p)/p)$.  Then
\begin{align}
\variation{K_{s,t}(a)-K_{s,t}(b)}{p}{T} &\leq \nonumber\\
2^{(p-1)/p}\bigg\{ \variation{\int_0^.(\nabla f(Y_{s,s+u}(a))&-\nabla f(Y_{s,s+u}(b)))\,K_{s,s+u}(a)\;dX_{s+u}}{p}{T}\nonumber\\
 + \variation{\int_0^. \nabla f(Y_{s,s+u}(b))(&K_{s,s+u}(a)-K_{s,s+u}(b))\;dX_{s+u}}{p}{T}\bigg\}\nonumber.
\end{align}
Let 
\begin{equation*}
h_t\dfn \int_0^t K_{s,s+u}(a)\;dX_{s+u},\quad g_t\dfn \int_0^t (K_{s,s+u}(a)-K_{s,s+u}(b))\;dX_{s+u}.
\end{equation*}
A bound of Young implies that
\begin{align}
\variation{h}{p}{T} &\leq  \big( 1+\zeta(1/p +1/q)\big) \variation{K_{s,s+.}(a)}{q}{T} \,\variation{X}{p}{T},\nonumber\\
\variation{g}{p}{T} &\leq  \big( 1+\zeta(1/p +1/q)\big) \variation{K_{s,s+.}(a)-K_{s,s+.}(b)}{q}{T} \,\variation{X}{p}{T}\nonumber.
\end{align}
So 
\begin{align}
\variation{K_{s,t}(a)-K_{s,t}(b)}{p}{T} &\leq\nonumber\\
2^{(p-1)/p}\variation{X}{p}{T}&\big(1+\zeta(1/p+1/q)\big)\nonumber\\
\times\;\bigg\{ \variation{\nabla f(Y_{s,s+u}(a))&-\nabla f(Y_{s,s+u}(b))}{r}{T}\variation{K_{s,s+.}(a)}{q}{T}\nonumber\\
 + \variation{\nabla &f(Y_{s,s+u}(b))}{r}{T}\variation{K_{s,s+.}(a)-K_{s,s+.}(b)}{q}{T}\bigg\}\nonumber.
\end{align}

Noting that 
\begin{align*}
&\variation{\nabla f(Y_{s,s+u}(a))-\nabla f(Y_{s,s+u}(b))}{q}{T}\\
\leq 2 &\variation{\nabla f(Y_{s,s+u}(a))-\nabla f(Y_{s,s+u}(b))}{q*}{T}^{q*/q}\norm{f}{\lip{\alpha}}^{(q-q*)/q} \norm{Y(a)-Y(b)}{\infty}^{(q-q*)/q}
\end{align*}
where $q*\in(p,q)$, we see that for $T$ sufficiently small 
$$
\variation{K_{s,t}(a)-K_{s,t}(b)}{p}{T}
$$ 
varies continuously in space.

We have not yet considered the inverse maps and their regularity.  The above proof can be modified to recover the result by considering the solution of the differential equation driven by the time-reversed path $(X_{t-u})_{0\leq u\leq t}$. 
\endproof

We note that the pathwise approach improves the following result of Kunita \cite{kun9} whenever the continuous integrator almost surely has finite $p$-variation for some $p<2$:

\begin{thm}[\cite{kun9} Theorem 2.3]
Let $f$ be a $C^k$ map with a bounded derivative, let $(x_t)_{t\geq0}$ be a square integrable martingale.  Then the solution of \eqref{e.1} defines a $C^{k-1}$ flow of diffeomorphisms for all t with probability one.
\end{thm}

The improvement is due to the increased smoothness of the integrator.  It is the pathwise approach that picks up the $p$-variation $(p<2)$ of the integrator thereby reducing the smoothness requirement on $f$.

The following example contains an interesting application of the above result.  The example, fractional Brownian motion, is a process which is not a semi-martingale \cite[Corollary 2.2]{lin}.  The process had been considered suitable for financial modeling due its heavy tail distribution.  (However Rogers \cite{rogers} has shown that arbitrage opportunities appear when using the fBm.)

\begin{eg}
Let $(X_t)_{t\geq0}$ be a fractional Brownian motion with Hurst index $H\in(\half,1)$.  That is, $X_t$ is a mean zero, Gaussian process whose covariance function is given by
\begin{align}
\Gamma^H(s,t) &= \half (t^{2H} + s^{2H} - \snorm{t-s}^{2H}) \, Var(X_1)\\
\intertext{where}
Var(X_1) &= {{1}\over{2H}} + \int_0^{\infty} \big( (1+s)^{H-\half}- s^{H-\half}\big)^2\;ds.
\end{align}
    
\begin{lem}\label{lem.fbm}
Fractional Brownian motion with Hurst index $H$ has finite $p$-variation almost surely for $p>1/H$.
\end{lem}

\proof
We estimate $\snorm{X_t - X_s}^p$ for $s<t\in\R^+$ as follows:
Partition the interval $[s,t]$ by selecting the largest dyadic interval $[(k-1)2^{-n}T,k2^{-n}]$ which is contained within $[s,t]$.  Then we add dyadic intervals to either side of the initial interval, which are chosen maximally with respect to inclusion in the interval $[s,t]$.  Continuing in this fashion we label the partition according to the lengths of the dyadics.  We note that there are at most two dyadics of the same length in the partition we can label these $[l_{1,k},r_{1,k}]$ and $[l_{2,k},r_{2,k}]$ where $r_{1,k}\leq l_{2,k}\,.$  Then
\begin{equation}
[u,v] = \bigcup_{k=1}^{\infty} \bigcup_{i=1,2} [l_{i,k},r_{i,k}].
\end{equation}
We use the H\"older inequality to show that
\begin{align}
\snorm{X_t-X_s}^p &\leq \bigg(\sum_{k=1}^{\infty} \sum_{i=1,2} \snorm{X_{l_{i,k}}-X_{r_{i,k}}}\bigg)^{p}\\
&\leq  \bigg( \sum_{n=1}^{\infty} n^{-\gamma/(p-1)} \bigg)^{p-1}\;\sum_{n=1}^{\infty} n^{\gamma} \bigg( \sum_{i=1,2} \snorm{X_{l_{i,k}}-X_{r_{i,k}}} \bigg)^{p}\\
&\leq C(p,\gamma) \;\sum_{n=1}^{\infty} n^{\gamma} \sum_{i=1,2} \snorm{X_{l_{i,k}}-X_{r_{i,k}}}^{p}.
\end{align}
One can uniformly bound $\snorm{X_t-X_s}^{p}$ for any pair of times $s<t \in [0,T]$ by extending the above estimate over all dyadic intervals at level $n$, that is
\begin{equation}
\snorm{X_t-X_s}^p \leq C(p,\gamma) \;\sum_{n=1}^{\infty} n^{\gamma} \sum_{i=1}^{2^n} \snorm{X_{l_{i,k}}-X_{r_{i,k}}}^{p}.
\end{equation}
The $p$-variation of $X_t$ can be estimated by this bound.
\begin{equation}\label{e.pbd}
\sup_{\pi} \sum_{\pi} \snorm{X_{t_{k+1}}-X_{t_k}}^p \leq C(p,\gamma) \;\sum_{n=1}^{\infty} n^{\gamma} \sum_{i=1}^{2^n} \snorm{X_{l_{i,k}}-X_{r_{i,k}}}^{p}.
\end{equation}
Noting that 
\begin{equation}
\expn{(X_t-X_s)^2} = \snorm{t-s}^{2H} \;Var(X_1).
\end{equation}
we see that the right hand side of \eqref{e.pbd} is finite a.s. whenever $Hp>1$.
\endproof

From Theorem \ref{thm.main} and Lemma \ref{lem.fbm} we deduce the following result:
\begin{cor}
Let $(X_t)_{t\geq0}$ be a fractional Brownian motion with Hurst index $H$.  Let $f$ be an $\alpha$-Lipschitz function for some $\alpha>1/H$.  Then the solution of \eqref{e.flo} driven by $(X_t)$ defines a stochastic flow of local $C^1$-diffeomorphisms. 
\end{cor} 
\end{eg}

\section{Discontinuous paths}\label{s.jump}
In this section we extend the results of Section \ref{s.flow} to discontinuous integrators.  The theory is applied to sample paths of certain classes of L\'evy processes, namely those that have finite $p$-variation a.s. for some $p<2$, see Section \ref{s.levy}.

First we determine the solution's behaviour when the integrator jumps.  There are two possibilities which we consider:  the first is an extension of the Lebesgue-Stieltjes integral, the second is based on a geometric approach.

Suppose that the discontinuous integrator has bounded variation.  Then the solution $Y$ would have a jump of the form
$$
Y_t-Y_{t-} = f(Y_{t-})\;(X_t-X_{t-})
$$
at a jump time of $X$.  If $X$ has finite $p$-variation for some $1<p<2$ we could insert the (countable) set of jumps of the above form at the discontinuities of $X$.  The disadvantage of this approach is that one cannot expect the solution to generate a flow of diffeomorphisms for a general $f$ \cite{m},\cite{lea2}.    

Another jump behaviour which we use is the following:  When a jump of the integrator occurs we insert fictitious time, allowing the jump to be traversed, thereby creating a continuous path.  Then we solve the differential equation driven by the continuous path.  Finally, we remove the fictitious time components of the integrator and the solution path.  We call this a geometric solution.  This jump behaviour has been considered before by \cite{marc} and \cite{kur2}.  The reason that we work with this second type of solution is that it experiences an 'instantaneous flow' along an integral curve at jump times.  This removes the flow problems associated to the discontinuous integrator \cite{m},\cite{lea2}.

Now we define a parametrisation for a c\`adl\`ag path $X$ of finite $p$-variation.  The parametrisation adds fictitious time, allowing the traversal of the discontinuities of $X$.  We prove that the resulting continuous path has the same $p$-variation as the original path.  We solve \eqref{e.1} driven by the continuous path from which we extract a geometric solution by removing the fictitious time (or undoing the parametrisation).
 
\begin{defn}\label{d.tau}
Let $X$ be a c\`adl\`ag path of finite $p$-variation.  Let $\delta>0\,,$ for each $n \geq 1\,,$ let $t_n$ be the time of the $n$'th largest jump of $X\,.$  We define a map $\tau^{\delta} : [0,T]\rightarrow [0, T+\delta
\sum_{i=1}^{\infty}\vert J(t_i) \vert^p] $ (where $J(u)$ denotes
the jump of the path $X$ at time u) in the following way:
\begin{equation}\label{e.tau}
\tau^{\delta}(t) = t + \delta \sum_{n=1}^{\infty} \vert J(t_n)\vert^p \chi_{\{t_n \leq t\}}(t).
\end{equation}
\end{defn}
The map $\tau^{\delta} : [0,T]\rightarrow[0,\tau^{\delta}(T)]$ extends the original time frame into one on which we define the continuous process $X^{\delta}(s)$ below.
\begin{defn}\label{d.xdel}
\begin{align}\label{e.xdel}
{x}^{\delta}(s)&\nonumber\\
=& \left\{
\begin{array}{ll}
 {x}(t) &\text{if}\; s = \tau^{\delta}(t) ,\\
 {x}(t_n^-) + (s - \tau^{\delta}(t_n^-)) j(t_n)\delta^{-1}\vert j(t_n) \vert ^{-p} &\text{if}\; s\in[\tau^{\delta}(t_n^-) \tau^{\delta}(t_n)) .
\end{array}
\right.
\end{align}
\end{defn}
\begin{rems}
We see that $(s,\,{X}^{\delta}_s),\;s\in [0,\tau^{\delta}(T)]$ is a parametrisation of the driving path ${X}\,.$  The terms $\vert J(t_n) \vert ^{p}$ in \eqref{e.tau} ensure that $\tau^{\delta}(t)$ remains finite.
\end{rems}
The following two results are proved in \cite{drew2}:

\begin{propn}\cite[Proposition 1.1]{drew2}\label{p.pvarn}
Let ${X}$ be a c\`adl\`ag path of finite $p$-variation.  Define a parametrisation of ${{X}}$ as above, then
\begin{eqnarray} \Vert {X}^{\delta} \Vert_{p,[0,\tau^{\delta}(T)]} =
\variation{{X}}{p}{T} \qquad \forall \delta > 0 .
\end{eqnarray}
\end{propn}
\begin{thm}\cite[Theorem 1.2]{drew2}\label{t.geo}
Let ${X}$ be a c\`adl\`ag path with finite $p$-variation for some $p<2\,.$  Let $f$ be a $\lip(\gamma)$ vector field on $\R^n$ for some $\gamma> p\,.$ Then there exists a unique geometric solution ${Y}$ with finite $p$-variation which solves the differential equation
\begin{equation}\label{eq.21}
d{Y}_t = f({Y}_t)\; d{X}_t \qquad {Y}_0 = {a} \in \R^n.
\end{equation}
\end{thm} 
We use Theorem \ref{thm.main} to deduce
\begin{cor}\label{cor.geoflo}
Let ${X}$ be a c\`adl\`ag path with finite $p$-variation for some $p<2\,.$  Let $f$ be a $\lip(\gamma)$ vector field on $\R^n$ for some $\gamma> p\,.$  Then the unique geometric solution to the following system of differential equations is a flow of local $C^1$-diffeomorphisms on $\R^d$:
\begin{equation}\label{e.geoflo}
Y_{s,t}(x) = Y_{s,s}(x) + \int_0^{t-s} f(Y_{s,s+u}(x))\;dX_{s+u}, \qquad Y_{s,s}(x)=x\in\R^d.
\end{equation}
\end{cor}

\section{$p$-variation of L{\'e}vy processes}\label{s.levy}

In this section we recall some results concerning the $p$-variation of L\'evy processes.  These determine when Corollary \ref{cor.geoflo} can be applied to sample paths of L\'evy processes.  

L\'evy processes are a combination of drift, diffusion and jumps.  The drift component has bounded variation whilst the Gaussian part has finite $p$-variation a.s. only for $p>2$.  Therefore we restrict our attention to the compensated, pure jump process given by:
\begin{equation}
\lim_{\epsilon\rightarrow0}\bigg\{\,\int_{\epsilon<\snorm{x}\leq 1} \bigg( x\;N_t(dx) - t x \; \nu(dx)\bigg) \,\bigg\}.
\end{equation}
To characterise the $p$-variation of L\'evy processes it is natural to look at the jump behaviour.  Blumenthal and Getoor measured this in the following way:
\begin{defn}\cite{blu2}
The index of a L\'evy process $X$ is defined to be 
\begin{equation}\label{e.index}
\beta \dfn \inf\bigg\{ \alpha >0 \,:\, \int_{\vert {y}\vert \leq 1} \vert {y}\vert^{\alpha}\, \nu(d{y}) < \infty\bigg\}.
\end{equation}
\end{defn}
It is a measure of the frequency of the process' small jumps.  If the integral
\begin{equation}\label{e.p1}
\int_{\snorm{x}\leq1} \snorm{x}\;\nu(dx) 
\end{equation} 
is finite then we see that the compensated jump process has bounded variation.  
  In \cite{blu1}, Blumenthal and Getoor proved the following result for the symmetric $\alpha$-stable processes $\alpha \in (0,2]\,:$
\begin{thm}\cite[Theorem 4.1]{blu1}
Let $\{X_t;\,t\geq0\}$ be the symmetric stable process in $\R^n$ of index $\alpha\in(0,2]\,.$  Then 
\begin{equation}
\prob{\variation{X}{\beta}{T} = \infty} = 1 \text{ or } 0
\end{equation}
depending on whether $\beta \leq \alpha$ or $\beta > \alpha\,.$
\end{thm}
Their proof centred on the special form that the characteristic functions of $\alpha$-stable L\'evy processes $(\alpha\in(0,2])$ take:
\begin{equation}
\expn{e^{iuX_t}} = e^{-\half t u^{\alpha}}
\end{equation}

It was conjectured that any pure jump L\'evy process would have finite $p$-variation provided $p>\beta$.  Br\'etagnolle \cite{bret} resolved the issue when the index $\beta<2$ with the following result:

\begin{thm}\label{bretthm} \cite{bret}
Let ${X}(t)$ be a process with independent increments with L\'evy measure $\nu\,.$  Then ${X}$ has finite $p$-variation for some $p<2$ on $[0,1]$ a.s. if and only if ${X}$ has no Gaussian part and
\begin{equation}\label{e.pbret}
\int _{\vert {x}\vert\leq 1} \vert {x}\vert ^p \, \nu (d{x})\, <\, \infty. 
\end{equation}
\end{thm}
\begin{rem}
In \cite{mon} the characterisation of the sample path $p$-variation of all L\'evy processes is completed.
\end{rem}

Using Theorem \ref{bretthm} we apply Corollary \ref{cor.geoflo} to a class of L\'evy processes.
\begin{cor}\label{cor.levy}
Let $(X_t)$ be a L\'evy process without a Gaussian part and with index $\beta$ less than two.  Let $f$ be an $\alpha$-Lipschitz vector field for some $\alpha>\beta$.  Then there exists a unique stochastic flow of local $C^1$-diffeomorphisms $(Y_{s,t})$ which solve the equation:
\begin{equation}
Y_{s,t}(x) = Y_{s,s}(x) + \int_0^{t-s} f(Y_{s,s+u}(x))\;dX_{s+u}, \qquad Y_{s,s}(x)=x.
\end{equation} 
\end{cor}

\noindent{\bf Acknowledgements}
The author wishes to thank his PhD. supervisor Terry Lyons. 

\bibliographystyle{plain}
\bibliography{institut,journal,publish,ref}
\address
\end{document}